\documentclass[12pt]{article}
 \usepackage{amsmath,amssymb,amsthm}
\usepackage{mathrsfs}
\usepackage{pgf,tikz}

 \def\draw #1 by #2 (#3){
  \vbox to #2{
    \hrule width #1 height 0pt depth 0pt
    \vfill
    \special{picture #3} 
    }
  }

 \def\scaleddraw #1 by #2 (#3 scaled #4){{
  \dimen0=#1 \dimen1=#2
  \divide\dimen0 by 1000 \multiply\dimen0 by #4
  \divide\dimen1 by 1000 \multiply\dimen1 by #4
  \draw \dimen0 by \dimen1 (#3 scaled #4)}
  }

\newtheorem{theorem}{Theorem}[section]
\newtheorem{example}{Example}
\newtheorem{problem}[example]{Problem}
\newtheorem{defin}[theorem]{Definition}
\newtheorem{lemma}[theorem]{Lemma}

\newtheorem{nt}{Note}

\newenvironment{pf}{\medskip\noindent{Proof:  \hspace*{-.4cm}}
       \enspace}{\hfill \qed \newline \medskip}

 \newcommand{\singlespacing}{\let\CS=\@currsize\renewcommand{\baselinestretch}{1}\tiny\CS}
 \newcommand{\oneandahalfspacing}{\let\CS=\@currsize\renewcommand{\baselinestretch}{1.25}\tiny\CS}
 \newcommand{\doublespacing}{\let\CS=\@currsize\renewcommand{\baselinestretch}{1.35}\tiny\CS}

 \newtheorem{rule-def}[theorem]{Rule}

\begin{document}
\baselineskip 16pt
 \newcommand{\la}{\lambda}
 \newcommand{\si}{\sigma}
 \newcommand{\ol}{1-\lambda}
 \newcommand{\be}{\begin{equation}}
 \newcommand{\ee}{\end{equation}}
 \newcommand{\bea}{\begin{eqnarray}}
 \newcommand{\eea}{\end{eqnarray}}

 \begin{center}
 {\Large \bf The minimum ABS index of trees with given number of
pendent vertices }\\
\end{center}

\begin{center}
{\bf V. Maitreyi$^a$, Suresh Elumalai$^{b}$, Selvaraj Balachandran$^c$}\\
$^a${\it Department of Mathematics, College of Engineering and Technology, Faculty of Engineering and Technology, 
SRM Institute of Science and Technology, \\ Kattankulathur, Chengalpet 603 203, India
}\\
$^b${\it Department of Mathematics, School of Humanities and Sciences, \\
SASTRA Deemed University, Thanjavur, India}
\end{center}

\baselineskip=0.30in
  \begin{center}
 {\Large \bf Abstract }\\
\end{center}
The recently developed atom-bond sum-connection (ABS) index is a variation of the well-studied graph-based molecular descriptors connectivity (Randi{\'c}), atom-bond connectivity, and sum-connectivity indices, defined by $$ABS(G) = \sum\limits_{uv \in E(G)} \sqrt{\frac{d_u + d_v - 2}{d_u + d_v}}.$$ We present the minimum $ABS$ index of trees with a given number of pendent vertices in this paper,  which provides a solution to  the  problem proposed by A. Ali, I. Gutman and  I. Redzepovi{\'c},  Atom-Bond Sum-Connectivity Index of Unicyclic Graphs and Some Applications, Electron. J. Math.(2023) 1-7.

\section{Introduction}
\indent Topological indices are numerical representations of chemical structures that are derived from their hydrogen-depleted graphs as a method for concise and efficient description of structural formulas that are used to research and forecast the relationships between structure and properties of organic compounds.  In this study, it is assumed that every graph is finite. We use common terminology and notation from graph theory.
Those that aren't defined here can be found in the books \cite{BM,CLZ,WW}.

\indent In 1975, Randi{\'c} \cite{MR} proposed  the following index
$$R(G) \sum\limits_{uv \in E(G)} \frac{1}{\sqrt{d(u) d(v)}}.$$
where the degree of vertex $u \in G $ is denoted by $d(u).$  
In structure-property and structure-activity relationship studies, the Randi{\'c} index is one of the most successful and effective molecular descriptors \cite{TC,DB,KH}. This descriptor's mathematical properties have been extensively studied. (see \cite{GF,LG,LS,LJ}).\\
\indent
Zhou and Trinajsti{\'c} \cite{ZT,BT} presented the sum-connectivity index $\chi(G)$  and the general sum-connectivity index $\chi_\alpha (G)$ of a graph $G$ and described them as follows:
 
$$\chi(G) = \sum\limits_{uv \in E(G)} \frac{1}{\sqrt{d(u) + d(v)}}$$
and $$\chi_\alpha (G) = \sum\limits_{uv \in E(G)} (d(u) + d(v))^\alpha$$
where $\alpha$ is a real number. It has been proved that the Randi{\'c} and sum-connectivity indices correlate closely with one another and with the $\pi$-electronic energy of benzenoid hydrocarbons\cite{LTZ}.

  Another variation of the Randi{\'c} index is the harmonic index of a graph. This index was introduced in \cite{SF} and is described as
$$H(G) = \sum\limits_{uv \in E(G)} \frac{2}{d(u) + d(v)}.$$  
More results and mathematical properties of harmonic index were discovered in the survey \cite{AZG}. In addition to the degree of the edge's end-vertices, the connectivity index was modified by taking into account the edge's degree $d(u) + d(v) - 2.$
The subsequent index was called atom-bond connectivity index (ABC) \cite{EES},  defined by
  $$ABC(G) = \sum\limits_{uv \in E(G)} \sqrt{\frac{d_G (u) + d_G (v) - 2}{d_G (u) d_G (v)}}.$$
\indent  For more information on the $ABC$ index and its chemical applications, one can refer \cite{ETRG,EES,GTRM,RF} and for  
mathematical properties, see the most recent review  \cite{ADDF}.\\

\indent Another invariant of $ABC$ index was put forward by Akbar Ali et al.., For a (molecular)graph $G$ its atom-bond sum-connectivity index is defined as $$ABS(G) = \sum\limits_{uv \in E(G)} \sqrt{\frac{d_u + d_v - 2}{d_u + d_v}} = \sum\limits_{uv \in E(G)} \sqrt{1- \frac{2}{d_u + d_v}}$$. 

In the study of topological indices, it is often interesting to investigate the extremal values and characterizing their corresponding graphs of any particular index  under certain parameters. Very recently, A. Ali, I. Gutman and  I. Redzepovi{\'c} proposed a problem regarding the minimum $ABS$ index of trees with a given number of pendent vertices \cite{AGR}. In this  paper, we give the extremal values and characterize the  corresponding extremal graphs for the minimum $ABS$ index of trees with a given number of pendent vertices.

\section{Preliminaries}

A modified version of the atom-bond connectivity index, utilizing the core idea of the sum-connectivity index and named atom-bond sum-connectivity (ABS) index, was recently put forward. The ABS index is defined as, $$ABS(G) = \sum\limits_{uv \in E(G)} \sqrt{\frac{d_u + d_v - 2}{d_u + d_v}} = \sum\limits_{uv \in E(G)} \sqrt{1- \frac{2}{d_u + d_v}}$$.
\indent Let $G$ be a graph. For any vertex $v \in V(G)$, we use $N_G (v)$ (or $N(v)$ if there is no ambiguity) to denote the set of neighbors of $v$ in $G$. Let $n_i$ be the number of vertices of degree $i$ in $G$. We use $E_2 (G)$ to denote the set of edges in $G$ such that for any edge $uv \in E_2 (G)$, we have $d(u) = d(v) = 2$. The maximum degree of $G$ is denoted by $\Delta(G)$. A pendent vertex is a vertex of degree 1 and a non-pendent vertex is a vertex of degree at least 2. An edge incident with a pendent vertex is called a pendent edge. We say a path $P = v_0v_1...v_s$ of $G$ is a pendent path if $v_0$ is a pendent vertex, $d(v_s) \geq 3$ and $d(v_1) =...= d(v_{s-1}) = 2$ (unless $s = 1$). We use $\mathcal{P}(G)$ to denote the set of all pendent paths in $G$. Let $G - uv$ be the graph obtained from $G$ by deleting the edge $uv \in E(G)$. We use $P_n$ and $S_n$ to denote the path and the star on $n$ vertices, respectively. We write $A := B$ to rename $B$ as $A$.\\ 
\indent A chemical tree is a tree with maximum degree at most 4. Let $\mathcal{T}_{n,k}$ and $\mathcal{CT}_{n,k}$ be the set of trees and chemical trees with
$n$ vertices and $k$ pendent vertices, respectively, where $2 \leq k \leq n - 1$. A tree $T$ is called a $(k, 3)$-regular tree if $T$ contains $k$ pendent vertices and every non-pendent vertex has degree 3 in $T$. It is easy to calculate that every $(k, 3)$-regular tree has $2k - 2$ vertices. We use $\mathcal{T}^*_{n,k}$ to denote the set of trees on $n$ vertices obtained from a $(k, 3)$-regular tree by adding at least one new vertex on each pendent edge. Note that each tree $T \in \mathcal{T}^*_{n,k}$ can also be stated as follows: $T$ is a tree on $n$ vertices, there are $k - 2$ vertices of maximum degree 3 in $T$ which induce a tree and any of these vertices is adjacent to either another vertex of degree 3 or a vertex of degree 2. Clearly, every tree in $\mathcal{T}^*_{n,k}$ contains $k $ pendent vertices, $n-2k+2$ vertices of degree 2 and $k - 2$ vertices of degree 3. Since $\mathcal{T}_{n,2} = \mathcal{CT}_{n,2} = \{ P_n \}, \mathcal{T}_{n,n-1}=\{ S_n\}$ and $\mathcal{CT}_{n,n-1} = \{ S_n | 3  \leq n \leq 5 \}$, we only consider $3 \leq k \leq n-2 $ in the following arguments. \\
\indent In order to prove our main results, we need the following three transformations on a tree $T$ which were introduced by Zhang, Lu and Tian. \\
\indent (i) Let $uv$ be an edge in $T$. If $T_{uv}$ is the graph obtained from $T$ by contracting the edge $uv$, i.e., identifying the two vertices $u$ and $v$ in $T- uv$, we say that $T_{uv}$ is obtained from $T$ by Transformation I. Hence $|V(T_{uv})| = |V(T)| - 1$. \\
\indent (ii) Let $v$ be a vertex in $T$ with $N(v) = V' \cup V''$ such that $V' \cap V'' = \phi, |V'| = s_1 \geq 1$ and $|V''| = s_2 \geq 1$. If $T_{v\rightarrow (s_1,s_2)}$ is the graph obtained from $T$ by splitting the vertex $v$ into two new vertices $v'$ and $v''$, connecting $v'$ and $v''$, and connecting $v'$ to all vertices in $V'$ and $v''$ to all vertices in $V''$, we say that $T_{v\rightarrow (s_1,s_2)}$ is obtained from $T$ by Transformation II. Clearly,
$|V(T_{v\rightarrow (s_1,s_2)})| = |V(T)| + 1$. \\ 
\indent (iii) Let $v$ be a vertex in $T$ with $d(v) = s \geq 4$. If $T_{v\rightarrow(3-reg)}$ is the graph obtained from $T$ by replacing the vertex $v$ with an
$(s, 3)$-regular tree $T'$ such that each vertex in $N(v)$ and each pendent vertex of $T'$ are identified one by one, we say that $T_{v\rightarrow (3-reg)}$ is obtained from $T$ by Transformation III (see Fig. 1 for an illustration). By the construction of $T_{v\rightarrow (3-reg)}$, it is easy to see that $|T_{v\rightarrow(3-reg)}| = |V(T)| + s - 3$.\\
We now give several auxiliary lemmas by applying the above transformations. These lemmas will then be used to prove our main results in Section 2.\\

\begin{tikzpicture}[x=1cm,y=1cm]
\draw [line width=0.5pt] (1.,9.)-- (2.,8.);
\draw [line width=0.5pt] (1.,8.)-- (2.,8.);
\draw [line width=0.5pt] (1.,7.)-- (2.,8.);
\draw [line width=0.5pt] (3.,8.)-- (2.,8.);
\draw [line width=0.5pt] (3.,8.)-- (2.5,7.);
\draw [line width=0.5pt] (3.,8.)-- (3.5,7.);
\draw [line width=0.5pt] (3.,8.)-- (4.,8.);
\draw [line width=0.5pt] (4.,8.)-- (5.,9.);
\draw [line width=0.5pt] (4.,8.)-- (5.,7.);
\draw [line width=0.5pt] (3.,8.)-- (3.,9.);
\draw [line width=0.5pt] (3.,9.)-- (3.,10.);
\draw [line width=0.5pt] (7.,8.)-- (8.,8.);
\draw [line width=0.5pt] (8.,8.)-- (9.,8.);
\draw [line width=0.5pt] (9.,8.)-- (10.,9.);
\draw [line width=0.5pt] (9.,8.)-- (10.,7.);
\draw [line width=0.5pt] (8.,9.)-- (8.,8.);
\draw [line width=0.5pt] (7.,10.)-- (8.,9.);
\draw [line width=0.5pt] (9.,10.)-- (8.,9.);
\draw [line width=0.5pt] (12.,9.)-- (13.,8.);
\draw [line width=0.5pt] (12.,8.)-- (13.,8.);
\draw [line width=0.5pt] (12.,7.)-- (13.,8.);
\draw [line width=0.5pt] (13.,8.)-- (14.,8.);
\draw [line width=0.5pt] (14.,9.)-- (14.,8.);
\draw [line width=0.5pt] (14.,8.)-- (15.,8.);
\draw [line width=0.5pt] (15.,8.)-- (16.,9.);
\draw [line width=0.5pt] (15.,8.)-- (16.,7.);
\draw [line width=0.5pt] (13.,10.)-- (14.,9.);
\draw [line width=0.5pt] (15.,10.)-- (14.,9.);
\draw [line width=0.5pt] (13.,11.)-- (13.,10.);
\draw [line width=0.5pt] (14.128554995509983,10.954140499741422)-- (15.,10.);
\draw [line width=0.5pt] (16.,11.)-- (15.,10.);
\draw (1.8837687977597792,8.00798368281827655) node[anchor=north west] {$v_1$};
\draw (3.017723212850804,9.0387612688057418) node[anchor=north west] {$v_2$};
\draw (3.2059757580802235,8.00740881022683288) node[anchor=north west] {$v$};
\draw (4.112982336818231985,8.338470208672717) node[anchor=north west] {$v_3$};
\draw (3.292206646796774,7.030635063138364) node[anchor=north west] {$v_4$};
\draw (2.2143205378398902,7.030635063138364) node[anchor=north west] {$v_5$};
\draw (6.770185825030985,8.884599170544204) node[anchor=north west] {$v'_1$};
\draw (7.230083898185923,10.422383352656027) node[anchor=north west] {$v'_2$};
\draw (9.112791635163946,10.465498797014302) node[anchor=north west] {$v'_3$};
\draw (10.133190484976463,9.473843576773968) node[anchor=north west] {$v'_4$};
\draw (10.176305929334738,7.461789506721118) node[anchor=north west] {$v'_5$};
\draw (12.863835294333903,8.827111911399838) node[anchor=north west] {$v^*_1$};
\draw (13.280617923130563,10.465498797014302) node[anchor=north west] {$v^*_2$};
\draw (15.278300178397323,10.479870611800393) node[anchor=north west] {$v^*_3$};
\draw (16.097493621204556,9.473843576773968) node[anchor=north west] {$v^*_4$};
\draw (16.154980880348923,7.418674062362842) node[anchor=north west] {$v^*_5$};
\draw (2.7891931292835617,5.938377139395388) node[anchor=north west] {$T$};
\draw (7.819328304415686,6.096467102042397) node[anchor=north west] {$T'$};
\draw (13.223130663986197,6.182697990758949) node[anchor=north west] {$T_{v\rightarrow(3-reg)}$};
\begin{scriptsize}
\draw [fill=black] (1.,9.) circle (2pt);
\draw [fill=black] (1.,8.) circle (2pt);
\draw [fill=black] (1.,7.) circle (2pt);
\draw [fill=black] (2.,8.) circle (2pt);
\draw [fill=black] (3.,8.) circle (2pt);
\draw [fill=black] (4.,8.) circle (2pt);
\draw [fill=black] (5.,9.) circle (2pt);
\draw [fill=black] (5.,7.) circle (2pt);
\draw [fill=black] (3.,9.) circle (2pt);
\draw [fill=black] (3.,10.) circle (2pt);
\draw [fill=black] (2.5,7.) circle (2pt);
\draw [fill=black] (3.5,7.) circle (2pt);
\draw [fill=black] (7.,8.) circle (2pt);
\draw [fill=black] (8.,8.) circle (2pt);
\draw [fill=black] (9.,8.) circle (2pt);
\draw [fill=black] (8.,9.) circle (2pt);
\draw [fill=black] (7.,10.) circle (2pt);
\draw [fill=black] (9.,10.) circle (2pt);
\draw [fill=black] (10.,9.) circle (2pt);
\draw [fill=black] (10.,7.) circle (2pt);
\draw [fill=black] (12.,9.) circle (2pt);
\draw [fill=black] (12.,8.) circle (2pt);
\draw [fill=black] (12.,7.) circle (2pt);
\draw [fill=black] (13.,8.) circle (2pt);
\draw [fill=black] (14.,8.) circle (2pt);
\draw [fill=black] (15.,8.) circle (2pt);
\draw [fill=black] (16.,9.) circle (2pt);
\draw [fill=black] (16.,7.) circle (2pt);
\draw [fill=black] (14.,9.) circle (2pt);
\draw [fill=black] (13.,10.) circle (2pt);
\draw [fill=black] (15.,10.) circle (2pt);
\draw [fill=black] (13.,11.) circle (2pt);
\draw [fill=black] (14.128554995509983,10.954140499741422) circle (2pt);
\draw [fill=black] (16.,11.) circle (2pt);
\end{scriptsize}
\end{tikzpicture}
\begin{center}
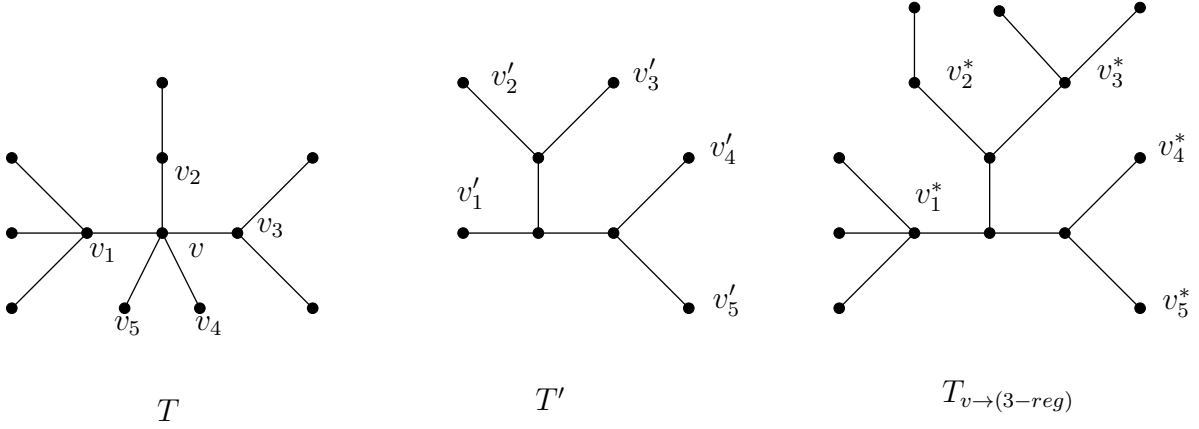
Figure 1. An Illustration of Transformation III.
\end{center}

\begin{lemma}
For any fixed $k \geq 1$, the function $f(x) = \sqrt{1-\frac{2}{x+k}} - \sqrt{1-\frac{2}{x}}$ is a strictly decreasing function for $x \geq 3$. \\
\begin{pf}
Since $k \geq 1$ and we have $$f'(x) = \frac{1}{(x+k-2)^{\frac{1}{2}} (x+k)^{\frac{3}{2}}} - \frac{1}{x^{\frac{3}{2}} (x-2)^{\frac{1}{2}}} < 0$$ for $x\geq 3$. This proves the lemma.
\end{pf}
\end{lemma}
\begin{lemma}
Let $T\in \mathcal{T} _{n,k}$ and let $v\in V(T)$ be a vertex of degree 2 such that the two neighbors of $v$ have degree at least 2 in $T$, then there exists a tree $T'$, such that $ABS(T) \geq ABS(T')$. \\
\begin{pf}
Let $N(v) = {u, w}$ with $d(u) = p \geq 2$ and $d(w) = q \geq 2$. Suppose that $x$ is a pendent vertex and $y$ is the unique neighbor of $x$ in $T$ with $d(y) = r \geq 2$. Let $T_{uv}$ be obtained from $T$ by Transformation I and let $T'$ be the graph obtained from $T_{uv}$ by adding a new edge to the pendent vertex $x$. Then $T' \in \mathcal{T}_{n,k}$ and

\begin{align*}
ABS(T) - & ABS(T')  =   \sqrt{\frac{p}{p+2}} + \sqrt{\frac{q}{q+2}} + \sqrt{\frac{r-1}{r+1}} - \sqrt{\frac{p+q-2}{p+q}} - \sqrt{\frac{r}{r+2}} - \sqrt{\frac{1}{3}}\\
& =   \Biggl[ \sqrt{\frac{p}{p+2}} + \sqrt{\frac{q}{q+2}} - \sqrt{\frac{p+q-2}{p+q}} \Biggr] - \Biggl[ \sqrt{\frac{r}{r+2}} - \sqrt{\frac{r-1}{r+1}} \Biggr]  - \sqrt{\frac{1}{3}}
\end{align*}
Let $$g(p,q) = \sqrt{\frac{p}{p+2}} + \sqrt{\frac{q}{q+2}} - \sqrt{\frac{p+q-2}{p+q}} $$ with $p,q \geq2.$ We have, 
$$\frac{\partial g(p,q)}{\partial p} = \frac{(p+2)^{-\frac{3}{2}}}{p^{\frac{1}{2}}} - \frac{(p+q)^{-\frac{3}{2}}}{(p+q-2)^{\frac{1}{2}}} \geq 0 $$
and $$\frac{\partial g(p,q)}{\partial q} = \frac{(q+2)^{-\frac{3}{2}}}{q^{\frac{1}{2}}} - \frac{(p+q)^{-\frac{3}{2}}}{(p+q-2)^{\frac{1}{2}}} \geq 0 $$
This implies that $g(p, q)$ is increasing for both $p$ and $q$ with $p, q \geq 2$. Then by Lemma 2.1 (with $k = 1$), we have 
\begin{eqnarray*}
ABS(T) - ABS(T') & = & g(p,q) - f(r+1) - \sqrt{\frac{1}{3}} \\
& \geq & g(2,2) - f(3) - \sqrt{\frac{1}{3}} \\
& = & 0
\end{eqnarray*}
So the assertion of the lemma holds.
\end{pf}
\end{lemma}
\begin{lemma}
Let $T \in \mathcal{T}_{n,k}$ with $E_2(T)\neq \phi$ and let $uv$ be a pendent edge in $T$ such that $d(u) = r \geq 3$ and $d(v) = 1$, then there exists a tree $T' \in \mathcal{T}_{n,k}$ such that $ABS(T) > ABS(T')$.\\
\begin{pf}
Since $E_2(T) \neq \phi$, let $T'$ be the graph obtained from $T$ by contracting an edge in $E_2(T)$ and adding a new edge to the pendent vertex $v$. Hence $T' \in \mathcal{T}_{n,k}$. Now by Lemma 2.1 (with $k = 1$), we have 
\begin{eqnarray*}
ABS(T) - ABS(T') & = & \sqrt{\frac{1}{2}} + \sqrt{1-\frac{2}{r+1}} - \sqrt{1-\frac{2}{r+2}} - \sqrt{\frac{1}{3}} \\
& = & \Biggl[ \sqrt{\frac{1}{2}} - \sqrt{\frac{1}{3}} \Biggr] - \Biggl[ \sqrt{1-\frac{2}{r+2}} - \sqrt{1-\frac{2}{r+1}} \Biggr] \\
& = & f(3) - f(r+1)   >  0
\end{eqnarray*}
Since $r\geq 3$This completes the proof of the lemma.
\end{pf}
\end{lemma}

\begin{lemma}
Let $T \in \mathcal{T}_{n,k}$ with $E_2(T) \neq \phi$ and let $v \in V(T)$ be a vertex of degree 4 such that $N(v) = \{u_1, u_2, u_3, u_4\}$, $d(u_1) \leq d(u_2) \leq d(u_3) \leq 3$ and $d(u_4) \leq 5$, then there exists a tree $T' \in \mathcal{T}_{n,k}$ such that $ABS(T) > ABS(T')$.
\end{lemma} 
\begin{pf}
Since $E_2(T) \neq \phi$, let $T'$ be the graph obtained from $T$ by contracting an edge in $E_2(T)$ and let $T' := T^*_{v\rightarrow (3-reg)}$ be
obtained from $T^*$ by Transformation III. Then $T' \in \mathcal{T}_{n,k}$. Since $d(u_1 ) \le d(u_2 ) \le d(u_3 ) \le 3$
 and $d(u_4) \leq 5$, by applying
Lemma 2.1 (with $k = 1$), we have $ABS(T) - ABS(T')  <  0.$
\begin{align*}
ABS(T) - ABS(T') & =  \sqrt{1-\frac{1}{2}} + \sum\limits_{i=1} ^4 \Biggl[ \sqrt{1-\frac{1}{4 + d(v_i)}} - \sqrt{1-\frac{2}{3 + d(v_i)}} \Biggr] - \sqrt{1-\frac{2}{6}} \\
& \leq  \sqrt{\frac{1}{2}} + 3\sqrt{\frac{5}{7}} + \sqrt{\frac{7}{9}} - 3\sqrt{\frac{2}{3}} - \sqrt{\frac{3}{4}} - \sqrt{\frac{2}{3}} \\
& <  0
\end{align*} 
\end{pf}
\begin{lemma}
Let $T \in \mathcal{T}_{n,k}$ and let $v \in V(T)$ be a vertex of degree $s \geq 4$ such that $N(v) = \{u_1, u_2,..., u_s\}$, $d(u_1) \leq d(u_2) \leq ... \leq d(u_s)$ and $|E_2(T)| \geq s - 3$. If $T$ satisfies one of the following properties:\\
(i) $s \geq 5$ and $d(u_{s-1}) \leq 3$; or \\
(ii) $s \geq 12$ and $d(u_{s-1}) \leq 4$\\
Then there exists a tree $T' \in \mathcal{T}_{n,k}$ such that $ABS(T) > ABS(T')$.
\\
\begin{pf}
Since $|E_2(T)| \geq s - 3$, let $T^*$ be the graph obtained from $T$ by contracting $s - 3$ edges in $E_2(T)$ and let $T' := T^*_{v\rightarrow(3-reg)}$ be obtained from $T^*$ by Transformation III. Hence $T' \in \mathcal{T}_{n,k}$. Let $d(u_{s-1}) = r$. Since $s \geq 4$, by Lemma 2.1
(with $k = s - 3 \geq 1$), we see that
\begin{align}
ABS(T) - ABS(T') & =   (s-3) \sqrt{\frac{1}{2}}+ \sum\limits_{i=1} ^s \Biggl[ \sqrt{1-\frac{2}{d(u_i) + s}} - \sqrt{1-\frac{2}{d(u_i) + 3}} \Biggr] - (s-3) \sqrt{\frac{2}{3}} \nonumber \\
 & \geq   (s-3)\Biggl[ \sqrt{\frac{1}{2}} - \sqrt{\frac{2}{3}} \Biggr] + (s-1) \Biggl[ \sqrt{1-\frac{2}{r + s}} - \sqrt{1-\frac{2}{r + 3}} \Biggr] \nonumber\\
 & +  \sqrt{1-\frac{2}{d(u_s) + s}} - \sqrt{1-\frac{2}{d(u_s) + 3}} \nonumber \\
 & >  (s-3) \Biggl[\sqrt{\frac{1}{2}} - \sqrt{\frac{2}{3}} \Biggr] + (s-1) \Biggl[ \sqrt{1-\frac{2}{r + s}} - \sqrt{1-\frac{2}{r + 3}} \Biggr]
\end{align}

Let $h(r,s) =  (s-3) \Biggl[\sqrt{\frac{1}{2}} - \sqrt{\frac{2}{3}} \Biggr] + (s-1) \Biggl[ \sqrt{1-\frac{2}{r + s}} - \sqrt{1-\frac{2}{r + 3}} \Biggr]$ where $r \geq 1$ and $s \geq 4$.  It is easy to see that $h(r, s)$ is strictly decreasing for $r$ by applying Lemma 2.1 (with $k = s - 3 \geq 1$). We have,
$$\frac{\partial h (r,s)}{\partial s} = \Biggl[ \sqrt{\frac{1}{2}} - \sqrt{\frac{2}{3}} \Biggr] + \sqrt{1-\frac{2}{r+s}} - \sqrt{1-\frac{2}{r+3}} + \frac{s-1}{\sqrt{1-\frac{2}{r+s}} (r+s)^2}$$
and 
$$\frac{\partial^2 h (r,s)}{\partial s^2} = \frac{2}{\sqrt{1-\frac{2}{r+s}} (r+s)^2} - \frac{s-1}{(1-\frac{2}{r+s})^2 (r+s)^4}  - \frac{2(s-1)}{\sqrt{1-\frac{2}{r+s}} (r+s)^3} > 0$$

\noindent Hence $\frac{\partial h (r,s)}{\partial s}$ is strictly increasing for $s$.\\
First, suppose that $T$ satisfies (i), i.e., $r \leq 3$ and $s \geq 5$. Then $h(r, s) \geq h(3, s)$ and 
$$\frac{\partial h(3,s)}{\partial s} \geq \frac{\partial h(3,s)}{\partial s}\Biggr|_{s=5} > 0 $$
This implies that $h(3, s)$ is strictly increasing for $s \geq 5$. Then by (1), we have 
$$ABS(T) - ABS(T')  >  h(r,s)  \geq h(3,s) \geq h(3,5) \geq 0 $$ 
So (i) holds. \\
Now assume that $T$ satisfies (ii). i.e.., $r\leq 4$ and $s\geq 12.$\\
Hence $h(r,s)  \geq h(4,s)$ and 
$$\frac{\partial h(4,s)}{\partial s} \geq \frac{\partial h(4,s)}{\partial s}\Biggr|_{s=12} > 0 $$
So $h(4,s)$ is strictly increasing for $s\geq 12.$ Therefore by (1), we know that,
$$ABS(T) - ABS(T')  >  h(r,s)  \geq h(4,s) \geq h(4,12) \geq 0 $$ 
This proves (ii) and hence completes the proof of the lemma.
\end{pf}
\end{lemma}
\noindent The final lemma in this section was proved by Zhang, Lu and Tian.
\begin{lemma} \cite{ZLT}
Let $T \in \mathcal{T}_{n,k}$ with $3 \leq k \leq \left\lfloor \frac{n+2}{3} \right\rfloor$ and $E_2(T) \subseteq E(\mathcal{P}(T))$, then
$$|E_2 (T)| \geq n_4 + 2n_5 +...+ (\Delta (T) - 3) n_{\Delta(T)}. $$
\end{lemma}

\section{Main Results}
In this section we prove the main results of this paper.
\begin{theorem}
Let $T \in \mathcal{T}_{n,k}$. If $3 \leq k \leq \left\lfloor \frac{n+2}{3} \right\rfloor$, then
$$ABS(T) \geq k\Biggl[ \sqrt{\frac{1}{3}}+\sqrt{\frac{3}{5}}\Biggr] + (n-3k+2)\sqrt{\frac{1}{2}} + (k-3)\sqrt{\frac{2}{3}}$$
with equality if and only if $T \in \mathcal{T}_{n,k} ^*$.
\end{theorem}
\begin{pf}
Suppose $T^*$ is the tree with the minimum ABS index among all trees in $T_{n,k}$ with $3 \leq k \leq \left\lfloor \frac{n+2}{3} \right\rfloor$. Since $k \geq 3$, we have $\Delta(T^*) \geq 3$. \\
We may assume that\\
(1) every vertex of degree 2 in $T^*$ is on a pendent path.\\
For otherwise, let $v_0$ be a vertex of degree 2 in $T^*$ such that $v_0$ is not on any pendent path. Then the two neighbors of $v_0$ must have degree at least 2 in $T^*$. By applying Lemma 2.2, there exists a tree $T_1 \in \mathcal{T}_{n,k}$ such that $ABS(T^*) \geq ABS(T_1)$. It follows from the proof of Lemma 2.2 that the above equality holds if and only if at least one neighbor of $v_0$ has degree 2 and every pendent vertex is adjacent to a vertex of degree 2 in $T^*$. By repeating the transformation in the proof of Lemma 2.2, we can construct a set of trees $\{T_i|i \geq 0 \}$ (with $T_0 := T^*$) in $\mathcal{T}_{n,k}$ and a set of vertices $\{v_i|i \geq 0\}$ such that $v_i$ is a vertex of degree 2 in $T_i$, $v_i$ is not on any pendent path and $ABS(T_i) \leq ABS(T_{i+1})$ for each $i \geq 0$. Note that the number of vertices of degree 2 which are not on pendent paths in $T_{i+1}$ is one less than that in $T_i$. Hence the above transformation will terminate in finite steps, i.e., there exists an integer $p \geq 0$ such that every vertex of degree 2 in $T_{p+1}$ is on a pendent path. This implies that $v_p$ is the unique vertex of degree 2 in $T_p$ which is not on any pendent path, and hence the two neighbors of $v_p$ have degree at least 3 in $T_p$. Therefore we have $ABS(T_0) \geq ABS(T_1) \geq ... \geq ABS(T_p) > ABS(T_{p+1})$. But this contradicts to the assumption that $T^*$ is the tree with the minimum ABS index among all trees in $\mathcal{T}_{n,k}$ with $3 \leq k \leq \left\lfloor \frac{n+2}{3} \right\rfloor$. So (1) holds. \\
\indent By (1), we know that $E_2(T^*) \subseteq E(\mathcal{P}(T^*))$. In the following argument, we will show $\Delta(T^*) = 3$. Suppose to the contrary that $\Delta(T^*) \geq 4$. Then by Lemma 2.6, we have
$$|E_2 (T^*)| \geq n_4 + 2n_5 +...+ (\Delta (T^*) - 3) n_{\Delta(T^*)} \geq \Delta(T^*) - 3 \geq 1. $$
Let $u_0$ be a vertex in $T^*$ with $d(u_0) = \Delta(T^*) \geq 4$ and let $P := u_0u_1...u_t$ be a path in $T^*$ with $d(u_t) \geq 4$. We choose $P$ such that the length of $P$ is as large as possible. Then we see that $t \geq 1$; otherwise, there exists a tree in $T_{n,k}$ with smaller ABS index than $T^*$ by Lemmas 2.4 and 2.5(i), contradicting the choice of $T^*$. Note that $d(u_i) \geq 3$ for each $1 \leq i \leq t - 1$ if $t \geq 2$ by (1).\\
\indent Let $N^*(u_{t-1}) := N(u_{t-1})$ if $t = 1$; otherwise, let $N^*(u_{t-1}) := N(u_{t-1}) \backslash \{u_{t-2}\}$ if $t \geq 2$. Clearly, $u_t \in N^*(u_{t-1})$. We claim that \\
(2) every vertex in $N^*(u_{t-1})$ has degree at most 4 in $T^*$ and $d(u_t) = 4$. \\
\noindent Let $x$ be a vertex in $N^*(u_{t-1})$ such that $N(x) = \{x_1, x_2,..., x_s\}$ and $d(x_1) \leq d(x_2) \leq ... \leq d(x_s)$. Then by the choice of $P$, we have $d(x_{s-1}) \leq 3$. Moreover, we may assume that $x_s = u_{t-1}$ since $d(u_{t-1}) = \Delta(T^*) \geq 4 > d(x_{s-1})$ if $t = 1$ ($u_{t-1} = u_0$ in this case) and $d(u_{t-1}) \geq 3 \geq d(x_{s-1})$ if $t \geq 2$. By applying Lemma 2.5(i) and by the choice of $T^*$, we deduce that $d(x) \leq 4$. Since $u_t  \in N^*(u_{t-1})$ and $d(u_t) \geq 4$, we further have $d(u_t) = 4$. This proves (2).\\
\indent We now consider the vertex $u_{t-1}$. By (2) and by the definition of $N^*(u_{t-1})$, we see that all the vertices in $N(u_{t-1})$ except $u_{t-2}$ (if $t \geq 2$) have degree at most 4 in $T^*$. Then by Lemma 2.5(ii) and by the choice of $T^*$, we have $d(u_{t-1}) \leq 7$. On the other hand, it follows from the proof of (2) that all the vertices in $N(u_t)$ except $u_{t-1}$ have degree at most 3 in $T^*$. Since $d(u_t) = 4$ (by (2)), we must have $d(u_{t-1}) \geq 6$; otherwise, there exists a tree in $\mathcal{T}_{n,k}$ with smaller ABS index than $T^*$ by Lemma 2.4, contradicting the choice of $T^*$. Hence we know that $6 \leq d(u_{t-1}) \leq 7$.\\
\indent Let $d(u_{t-1}) = r$ and let $N(u_{t-1}) = \{w_1, w_2,..., w_r\}$ with $w_r := u_{t-2}$ if $t \geq 2$. Then by (2), $d(w_i) \leq 4$ for each $1 \leq i \leq r - 1$. Since $|E_2(T^*)| \geq 1$, let $\bar{T}$ be the graph obtained from $T^*$ by contracting an edge in $E_2(T^*)$ and let $T' := \bar{T}_{u_{t-1}\rightarrow(3,r-3)}$ be obtained from $\bar{T}$ by Transformation II. Then $T' \in \mathcal{T}_{n,k}$. Since $6 \leq r \leq 7$, by Lemma 2.1 (with $k = r - 4 > 1$ and $k = 2$,
respectively), we have
\begin{align}
 ABS(T^*) & - ABS(T')  =  \sqrt{1-\frac{2}{4}} + \sum\limits_{i=1}^3 \Biggl[ \sqrt{1-\frac{2}{d(w_i) + r}} - \sqrt{1-\frac{2}{d(w_i) + 4}} \Biggr] \nonumber \\
& +   \sum\limits_{i=4}^r \Biggl[ \sqrt{1-\frac{2}{d(w_i) + r}} - \sqrt{1-\frac{2}{d(w_i) + r - 2}} \Biggr] - \sqrt{1-\frac{2}{r+2}} \nonumber \\
& \geq   \sqrt{\frac{1}{2}} + 3 \Biggl[ \sqrt{1-\frac{2}{r + 4}} - \sqrt{1-\frac{2}{8}} \Biggr] + (r-4) \Biggl[ \sqrt{1-\frac{2}{r+4}} - \sqrt{1-\frac{2}{r+2}} \Biggr] \nonumber \\
& +   \Biggl[ \sqrt{1-\frac{2}{d(w_r) + r}} - \sqrt{1-\frac{2}{d(w_r) + r - 2}} \Biggr] - \sqrt{\frac{r}{r+2}} \nonumber \\
& >   \sqrt{\frac{1}{2}} + 3 \Biggl[ \sqrt{\frac{r+2}{r + 4}} - \sqrt{\frac{3}{4}} \Biggr] + (r-4) \Biggl[ \sqrt{\frac{r+2}{r+4}} - \sqrt{\frac{r}{r+2}} \Biggr] - \sqrt{\frac{r}{r+2}} \nonumber \\
&= \sqrt{\frac{1}{2}} + (r-1) \sqrt{\frac{r+2}{r + 4}} - 3 \sqrt{\frac{3}{4}} - (r-3) \sqrt{\frac{r}{r+2}}  
\end{align}

Let $$\psi (r) = \sqrt{\frac{1}{2}} + (r-1) \sqrt{\frac{r+2}{r + 4}} - 3 \sqrt{\frac{3}{4}} - (r-3) \sqrt{\frac{r}{r+2}} $$
If $r=7$, then $$\psi (r) = \sqrt{\frac{1}{2}} + 6 \sqrt{\frac{9}{11}} - 3 \sqrt{\frac{3}{4}} - 4 \sqrt{\frac{7}{9}} \geq 0$$
If $r=8$, then $$\psi (r) = \sqrt{\frac{1}{2}} + 7 \sqrt{\frac{5}{6}} - 3 \sqrt{\frac{3}{4}} - 4 \sqrt{\frac{4}{5}} \geq 0$$

\noindent Then by (2), in both cases, we have $$ABS(T^*) - ABS (T') > \psi (r) \geq 0.$$
\noindent But this contradicts to the assumption that $T^*$ is the tree with the minimum ABS index among all trees in $\mathcal{T}_{n,k}$. Hence we have \\
(3) $\Delta(T^*) = 3.$ \\
\indent Now by (3), we have $n_1 + n_2 + n_3 = n$ and $n_1 + 2n_2 + 3n_3 = 2(n - 1)$. Since $n_1 = k$ and $3 \leq k \leq \left\lfloor \frac{n+2}{3} \right\rfloor$, we know that $n_2 = n - 2k + 2 \geq k$. Suppose that there exists a pendent vertex in $T^*$ whose neighbor has degree 3. Then $E_2(T^*)\neq \phi$ since $n_2 \geq k$ and by (1). By Lemma 2.3 and by the choice of $T^*$, we derive a contradiction. So each pendent vertex in $T^*$ is adjacent to a vertex of degree 2, and hence each pendent path has length at least 2 in $T^*$. This implies that $T^* \in \mathcal{T}_{n,k}^*$. It is easy to check
that $$ABS(T^*) = k\Biggl[ \sqrt{\frac{1}{3}}+\sqrt{\frac{3}{5}}\Biggr] + (n-3k+2)\sqrt{\frac{1}{2}} + (k-3)\sqrt{\frac{2}{3}}$$
by a simple calculation. This completes the proof of the theorem.
\end{pf}
Since $T^*_{n,k} \subseteq \mathcal{CT}_{n,k} \subseteq T_{n,k}$, we immediately obtain the following result for chemical trees.
\begin{theorem}
Let $T \in \mathcal{CT}_{n,k}$ with $3 \leq k \leq \left\lfloor \frac{n+2}{3} \right\rfloor$. Then
$$ABS(T) \geq k\Biggl[ \sqrt{\frac{1}{3}}+\sqrt{\frac{3}{5}}\Biggr] + (n-3k+2)\sqrt{\frac{1}{2}} + (k-3)\sqrt{\frac{2}{3}}$$
\end{theorem}

\end{document}